\documentclass[twoside]{article}
\title{A generalization of  strongly Gorenstein projective modules}
\date{}
\usepackage{amsfonts}
\usepackage{amsmath}
\usepackage{amssymb}
\usepackage{latexsym}
\newtheorem{thm}{\bf Theorem}[section]
\newtheorem{cor}[thm]{\bf Corollary}

\newtheorem{prop}[thm]{\bf Proposition}
\newtheorem{defn}[thm]{\bf Definition}

\newtheorem{exmp}[thm]{\bf Example}

\catcode`\ç=13
\defç{\c{c}}
\catcode`\é=13
\defé{\'e}
\catcode`\à=13
\defà{\`a}
\catcode`\è=13
\defè{\`e}
\catcode`\â=13
\defâ{\^a}
\catcode`\ù=13
\defù{\`u}
\catcode`\ê=13
\defê{\^e}
\catcode`\î=13
\defî{\^\i}
\catcode`\ô=13
\defô{\^o}

\def\proof{{\parindent0pt {\bf Proof.\ }}}
\def\pd{{\rm pd}}
\def\Im{{\rm Im}}
\def\Ker{{\rm Ker}}
\def\Ext{{\rm Ext}}
\def\Tor{{\rm Tor}}
\def\Hom{{\rm Hom}}
\newcommand{\cqfd}
{\hspace{1cm}
\rule{2mm}{2mm}%
\medbreak%
\par%
}

\begin{document}
\thispagestyle{empty}
\maketitle \vspace*{-2cm}
\begin{center}{\large\bf Driss Bennis and Najib Mahdou}

\bigskip
\small{Department of Mathematics, Faculty of Science and
Technology of Fez,\\ Box 2202, University S. M.
Ben Abdellah Fez, Morocco, \\[0.12cm] driss\_bennis@hotmail.com\\
mahdou@hotmail.com}
\end{center}

\bigskip\bigskip
\noindent{\large\bf Abstract.} This paper generalize the idea of the
authors in J. Pure Appl. Algebra \textbf{210} (2007) 437--445.
Namely, we define and study a particular case of Gorenstein
projective modules. We investigate some change of rings results for
this new kind of modules. Examples over not necessarily Noetherian
rings are given.\bigskip

\footnotetext[0]{ *Journal of Algebra and Its Applications, To appear 
}\value{counter}\setcounter{footnote}{0}

\medskip

\begin{section}{Introduction}  Throughout this paper, $R$ denotes
a non-trivial associative ring with identity, and all modules are,
if not specified otherwise, left $R$-modules. For an $R$-module $M$,
we use $\pd_R(M)$  to denote the classical projective dimension  of
$M$. It is convenient to use ``local" to refer to (not necessarily
Noetherian) rings with a unique maximal ideal.  \bigskip

We assume that the reader is familiar with the Gorenstein
homological dimension theory. Some references are \cite{LW, CFH,
Rel-hom, HH}. Nevertheless,  it is convenient to give a brief
history of the Gorenstein dimension theory.\bigskip

In the sixties,  Auslander and Bridger  introduced an homological
dimension, called G-dimension, for finitely generated modules over
Noetherian rings \cite{A1, A2}. Several decades later, this
homological dimension was extended, by Enochs et al. \cite{GoInPj,
GoIn, Fox}, to Gorenstein projective dimension of modules that are
not necessarily finitely generated and over not necessarily
Noetherian rings. And, dually, they defined the Gorenstein
injective dimension. Then, to complete the analogy with the
classical homological dimensions, Enochs  et al. \cite{GoPlat}
introduced the Gorenstein flat dimension.\medskip

\indent In the last years, the Gorenstein homological dimensions
have become a vigorously active area of research (see \cite{LW} and
\cite{Rel-hom} for more details). In 2004, Holm \cite{HH}
generalized several results which are already obtained over
Noetherian rings to associative rings. Recently, in \cite{BM}, the
authors introduced a particular case of Gorenstein projective,
injective, and flat modules, which are defined, respectively, as
follows:

\begin{defn}[\cite{BM}]\label{defSG}
\begin{itemize}
\item \textnormal{ An $R$-module $M$ is said to be strongly
Gorenstein projective (SG-projective for short), if there exists an
exact sequence of  projective modules $$ \mathbf{P}=\
\cdots\stackrel{f}{\longrightarrow}P\stackrel{f}{\longrightarrow}P\stackrel{f}{\longrightarrow}P
\stackrel{f}{\longrightarrow}\cdots $$ such that  $M \cong \Im(f)$
and such that $\Hom_R ( -, Q) $ leaves the sequence $\mathbf{P}$
exact whenever $Q$ is a projective $R$-module.}
        \item  \textnormal{The strongly Gorenstein injective
(SG-injective for short) modules are defined dually.}
        \item \textnormal{
An $R$-module $M$ is said to be strongly Gorenstein flat (SG-flat
for short), if there exists an exact sequence of  flat $R$-modules
 $$
\mathbf{F}=\ \cdots\stackrel{f}{\longrightarrow}F
\stackrel{f}{\longrightarrow}F\stackrel{f}{\longrightarrow}F
\stackrel{f}{\longrightarrow}\cdots
    $$ such that  $M \cong \Im(f)$ and such that $I \otimes_R - $ leaves the sequence
$\mathbf{F}$ exact whenever $I$ is an injective right $R$-module.}
\end{itemize}
\end{defn}

It is proved that the class of all strongly Gorenstein projective
modules is an intermediate class between the ones of projective
modules and Gorenstein projective modules \cite[Proposition
2.3]{BM}; i.e., we have the following inclusions
\begin{eqnarray*}
\{projective\ modules\}  &\subseteq&  \{SG\!-\!projective\ modules\}\\
  &\subseteq&  \{G\!-\!projective\ modules\}
\end{eqnarray*}
which are  strict by \cite[Examples 2.5 and 2.13]{BM}. The
principal role of the strongly Gorenstein projective, injective,
and flat modules is to give a simple characterization of
Gorenstein projective, injective, and flat modules, respectively
(see \cite[Theorems 2.7 and 3.5]{BM}). For example, we have, in
the Gorenstein projective case,  that a module is Gorenstein
projective if and only if it is a direct summand of a strongly
Gorenstein projective module\,\footnote{\;  In  \cite{BM} the
rings are assumed to be commutative. However, one can show easily
that this assumption is not necessary for some results such as the
principal one  \cite[Theorem 2.7]{BM} and
  \cite[Proposition 2.9]{BM}, which we need in this paper.}. The strongly
Gorenstein modules confirm that there is an analogy between the
notion of \textquotedblleft Gorenstein projective, injective, and
flat modules\textquotedblright  and the notion of the usual
\textquotedblleft projective, injective, and flat
modules\textquotedblright. This is, in fact obtained, because the
strongly Gorenstein projective, injective, and flat modules have
simpler characterizations than their  Gorenstein correspondent
modules \cite[Propositions  2.9 and 3.6]{BM}. For instance, a module
$M$ is  strongly Gorenstein projective if and only if there exists a
short exact sequence $0\rightarrow M\rightarrow P\rightarrow
M\rightarrow 0$, where $P$ is a projective  module, such that the
short sequence   $0\rightarrow \Hom(M,Q)\rightarrow
\Hom(P,Q)\rightarrow \Hom(M,Q)\rightarrow 0$ is exact for any
projective module $Q$.\bigskip

The aim of this paper is to generalize the notion of
\textquotedblleft strongly Gorenstein projective
modules\textquotedblright to $n$-strongly Gorenstein projective
modules (Definition \ref{DefnSGproj}). We generalize some results
of \cite{BM} (see Theorem \ref{thm-car-n-SG-proj}  and Proposition
\ref{pro-SG-pdim}), we establish some change of rings results
(Theorem \ref{thm-ch-ring-SGp} and Corollary
\ref{cor-ch-ring-SGp}), and we give some examples of $n$-strongly
Gorenstein projective modules  over not necessarily Noetherian
rings (see Examples \ref{exm1} and \ref{exm2}).

\end{section}
\begin{section}{$n$-Strongly
Gorenstein projective modules} In this section, we introduce and
study $n$-strongly Gorenstein projective modules which are defined
as follows:

\begin{defn}\label{DefnSGproj}\textnormal{Let $n$ be a positive integer.
  An $R$-module $M$ is said to be $n$-strongly
Gorenstein projective ($n$-SG-projective for short), if there exists
an exact sequence of  $R$-modules $$0\rightarrow M\rightarrow
P_n\rightarrow\cdots\rightarrow P_1 \rightarrow M\rightarrow 0,$$
where each $P_i$ is projective, such that $\Hom_R ( -, Q) $ leaves
the sequence  exact whenever $Q$ is a projective $R$-module.}
\end{defn}

Consequently, the $1$-strongly Gorenstein projective   modules are
just the strongly Gorenstein projective  modules (by
\cite[Proposition 2.9]{BM}).\bigskip

Recall that a Gorenstein projective module is projective if and
only if it has finite projective dimension \cite[Proposition
2.27]{HH}. For the $n$-strongly Gorenstein projective modules, we
have the following result, which is a generalization of
\cite[Corollary 2.11]{BM}.

\begin{prop}\label{pro-SG-pdim} Let $n$ be a positive integer.
An $n$-strongly Gorenstein projective module is projective if and
only if it has finite flat dimension.
\end{prop}
\proof The ``only if'' part is trivial.\\
We prove the  ``if'' part. Let $M$ be an $n$-strongly Gorenstein
projective module with finite flat dimension. Then, there exists  an
exact sequence of modules $$(\alpha)\quad 0\rightarrow M\rightarrow
P_n\rightarrow\cdots\rightarrow P_1 \rightarrow M\rightarrow 0,$$
where each $P_i$ is  projective, and, by standard arguments, $M$ is
flat.  We decompose the above sequence $(\alpha)$ into short exact
sequences
$$(\alpha_i)\quad 0\rightarrow H_{i}\rightarrow      P_i\rightarrow
H_{i-1}\rightarrow 0,$$ where  $H_{n}= M=H_0$ and
$H_{i}=\Ker(P_{i}\rightarrow H_{i-1})$ for $i=n-1,...,1$. It is
evident that each $H_i$ is flat, so is the finite direct sum
$H=\bigoplus\limits_{1\leq i \leq n}H_i$. Now, by adding the short
exact sequences $(\alpha_i)$, we obtain the following short exact
sequence:
$$ 0\rightarrow H\rightarrow      P \rightarrow
H \rightarrow 0,$$ where $P=\bigoplus\limits_{1\leq i \leq n}P_i$.
Finally, \cite[Theorem 2.5]{Period} gives the desired result.\cqfd

\begin{cor}
If  $R$ has  finite weak dimension, then the class of all
projective $R$-modules and the class of all Gorenstein projective
$R$-modules are the same class.
\end{cor}
\proof Apply  Proposition \ref{pro-SG-pdim} and \cite[Theorem
2.7]{BM}.\cqfd\bigskip

The following gives an  example, over  a not necessarily
Noetherian ring,  of a  $1$-strongly Gorenstein projective module
which is not projective.

\begin{exmp}\label{exm1}
Let $R$ be a local ring and consider the ring $S=R[[X]]/(X^{2})$.
Denote by $\overline{X}$   the residue class in $S$ of $X$. Then,
the ideal $(\overline{X})$ is strongly Gorenstein projective but
it is not projective.
\end{exmp}
\proof  First, since $R$ is local,  $S$ is also local, and so the
ideal $(\overline{X})$ is not  projective (since  $(\overline{X})^2=0)$.\\
Now, we prove that $(\overline{X})$ is strongly Gorenstein
projective.\\
Let $x$ be the homothety given by multiplication by
$\overline{X}$. We have the exact sequence:
             $$ \cdots \longrightarrow S
             \stackrel{x}{\longrightarrow} S
             \stackrel{x}{\longrightarrow} S \longrightarrow\cdots
             $$ such that $\Ker\, x=\Im\,x=(\overline{X})$.\\
             Since $S$ is local, every projective
              $S$-module is free. Then, consider a free $S$-module
              $S^{(I)}$, where $I$ is an index set, and applying the functor
$\Hom_S(-,S^{(I)} )$ to the above sequence, we get the following
sequence: $$   \cdots  \longrightarrow  \Hom(S, S^{(I)})
\stackrel{\Hom_S(x,S^{(I)})}{\longrightarrow} \Hom(S,S^{(I)})
 \longrightarrow  \cdots$$
By the natural isomorphism $ \Hom(S, S^{(I)})
\stackrel{\cong}{\longrightarrow} S^{(I)}$, we get the following
commutative diagram:
 $$
    \begin{array}{ccccccc}
      \cdots & \longrightarrow & \Hom(S, S^{(I)}) &
\stackrel{\Hom_S(x,S^{(I)})}{\longrightarrow} & \Hom(S,S^{(I)})
&\longrightarrow& \cdots \\
        &   &\cong\downarrow &  & \cong\downarrow  &  & \\
      \cdots &   \longrightarrow &S^{(I)} &
  \stackrel{x}{\longrightarrow} &S^{(I)}&   \longrightarrow& \cdots
    \end{array}
    $$
Evidently, the lower sequence in the diagram above is exact, then
so is the upper sequence, which means that the ideal
$(\overline{X})$ is strongly Gorenstein projective.\cqfd

\begin{prop}\label{Pro-relation} Let $n$ be a positive integer.
\begin{enumerate} \item Every strongly Gorenstein projective
  module is $n$-strongly Gorenstein projective.
        \item  Every $n$-strongly Gorenstein projective
 module  is  Gorenstein projective.
\end{enumerate}
\end{prop}
\proof 1. Let $M$ be a  strongly Gorenstein  projective
 module. Then, there exists, by   \cite[Proposition  2.9]{BM},
a short exact sequence:  $0\longrightarrow M
\stackrel{f}\longrightarrow P \stackrel{g}\longrightarrow
M\longrightarrow 0$, where $P$ is a projective module,  such that
the short sequence $0\rightarrow \Hom(M,Q)\rightarrow
\Hom(P,Q)\rightarrow \Hom(M,Q)\rightarrow 0$ is exact for any
projective module $Q$. Thus, we get an exact sequence:
$$0\longrightarrow M \stackrel{f}\longrightarrow P\stackrel{f\circ g}\longrightarrow
\cdots\stackrel{f\circ g}\longrightarrow
 P \stackrel{g}\longrightarrow M\longrightarrow 0,$$
such that $\Hom ( -, Q) $ leaves the sequence exact whenever $Q$
is a projective module. This, shows that $M$ is  $n$-strongly
Gorenstein projective for every positive integer $n$.\\
2. Let $M$ be a   $n$-strongly Gorenstein projective module. Then,
there exists an exact sequence $$0\longrightarrow M
\stackrel{f_{n+1}}\longrightarrow
P_n\stackrel{f_{n}}\longrightarrow\cdots\stackrel{f_2}\longrightarrow
P_1 \stackrel{f_1}\longrightarrow M\longrightarrow 0,$$ where each
$P_i$ is a projective module, such that $\Hom ( -, Q) $ leaves the
sequence exact whenever $Q$ is a projective module. Thus, we get
an exact sequence:
 $$\cdots\stackrel{f_2}\longrightarrow P_1\stackrel{f_{n+1}\circ f_1}\longrightarrow
P_n\stackrel{f_{n}}\longrightarrow\cdots\stackrel{f_2}\longrightarrow
P_1 \stackrel{f_{n+1} \circ f_1}\longrightarrow
P_n\stackrel{f_{n}}\longrightarrow\cdots, $$ such that $\Hom ( -,
Q) $ leaves the sequence exact whenever $Q$ is a projective
module. This shows that $M$ is Gorenstein projective.\cqfd\bigskip

The following is an extension of \cite[Example 2.13]{BM}, such
that we give an example of $2$-strongly Gorenstein projective
modules, over a not necessarily Noetherian ring, which are not
$1$-strongly Gorenstein projective.

\begin{exmp}\label{exm2}
Let $R$ be a local ring and consider the ring $S=R[[X ,Y]]/(XY)$.
Then:
\begin{enumerate}
    \item The two ideals $(\overline{X })$ and $(\overline{Y})$
are $2$-strongly Gorenstein projective, where $\overline{X }$ and
$\overline{Y }$ are the residue classes in $S$ of $X $ and $Y$,
respectively.
    \item $(\overline{X })$ and $(\overline{Y})$ are not strongly
Gorenstein projective.
\end{enumerate}
\end{exmp}
\proof The proof of  the  statement $(1)$ is similar to the one of
Example \ref{exm1}. We only need to note that we have two exact
sequences: $$
\begin{array}{ccccccccccc}
  0&\longrightarrow &(\overline{X})&\longrightarrow & R& \stackrel{y}\longrightarrow & R &
 \stackrel{x}\longrightarrow &  (\overline{X }) & \longrightarrow & 0 \\
0&\longrightarrow &(\overline{Y})&\longrightarrow & R &
\stackrel{x}\longrightarrow & R &  \stackrel{y} \longrightarrow &
(\overline{Y }) & \longrightarrow &  0
\end{array}$$
where $x$ and $y$ are the multiplications by $\overline{X}$ and
$\overline{Y}$, respectively.\\
The proof of the statement $(2)$ is the same as the one of
\cite[Example 2.13 (2)]{BM}.\cqfd\bigskip

Obviously, any example of an injective module which is not
projective is an example of (strongly) Gorenstein injective module
which is not (strongly) Gorenstein projective.  From the
following, which is a generalization of \cite[Exercise 4.2, p.
115]{Rot}, we can construct an example of an ($n$-strongly)
Gorenstein injective module which is neither injective nor
($n$-strongly) Gorenstein projective.

\begin{prop}\label{pro-domain}
If $R$ is an integral domain but not a field, an $R$-module that
is simultaneously Gorenstein projective and Gorenstein injective
must be 0.
\end{prop}
\proof Assume that an $R$-module $M$ is simultaneously Gorenstein
projective and Gorenstein injective. Then, as a Gorenstein
projective $R$-module, $M$ embeds in a projective $R$-module, and
so it is torsion-free (since $R$ is an integral domain). On the
other hand, as a Gorenstein injective  $R$-module, $M$ is a
quotient of an injective
 $R$-module, and so it is divisible (by  \cite[Exercise 3.14, p.
69]{Rot}). Thus, by  \cite[Exercise 3.19, p. 70]{Rot},  $M$ is an
injective  $R$-module, then it is also projective (since it embeds
in a projective $R$-module). Finally, by \cite[Exercise 4.2, p.
115]{Rot}, $M$ must be 0.\cqfd\bigskip

The next result is a characterization of $n$-strongly Gorenstein
projective modules. It is a generalization of \cite[Proposition
2.9]{BM}.

\begin{thm}\label{thm-car-n-SG-proj} For any module $M$  and any
positive integer $n$, the following are equivalent:
\begin{enumerate}
  \item $M$ is  $n$-strongly Gorenstein projective;
  \item There exists an exact sequence
$0\rightarrow M\rightarrow P_n\rightarrow\cdots\rightarrow P_1
\rightarrow M\rightarrow 0$, where each $P_i$ is a projective
module, such that $\Hom  ( -, Q') $ leaves the sequence  exact
whenever $Q'$  is a  module with finite projective dimension;
  \item There exists an exact sequence
$0\rightarrow M\rightarrow P_n\rightarrow\cdots\rightarrow P_1
\rightarrow M\rightarrow 0$, where  each $P_i$ is a projective
module, and there exists a positive integer $j$ such that
$\Ext^{j+1}(M,Q)=\Ext^{j+2}(M,Q)=\cdots=\Ext^{j+n}(M,Q)=0$  for any
projective module $Q$;
  \item There exists an exact sequence
$0\rightarrow M\rightarrow P_n\rightarrow\cdots\rightarrow P_1
\rightarrow M\rightarrow 0$, where each $P_i$ is a projective
module, and there exists a positive integer $j$ such that
$\Ext^{j+1}(M,Q')=\Ext^{j+2}(M,Q')=\cdots=\Ext^{j+n}(M,Q')=0$  for
any module $Q'$ with finite projective dimension.
\end{enumerate}
\end{thm}
\proof $1\Leftrightarrow 2 $. Obviously the second condition is
stronger than the first, then only the implication $1 \Rightarrow
2$ merits a proof. So, let $M$ be an $n$-strongly Gorenstein
projective module. Then, there exists an exact sequence
$$\mathbf{P}=\quad 0 \rightarrow M\rightarrow
P_n\rightarrow\cdots\rightarrow P_1 \rightarrow M\rightarrow 0,$$
where each $P_i$ is a projective module, such that $\Hom  ( -, Q)
$ leaves the sequence  exact whenever $Q$ is a projective module.
Let $Q'$ be a module with $\pd(Q')=m<\infty$ and consider a short
exact sequence $0\rightarrow K \rightarrow F\rightarrow
Q'\rightarrow 0$, where $F$ is a free module, then $\pd(K)=m-1$.
Thus,
$$0\rightarrow \Hom  ( \mathbf{P}, K)  \rightarrow \Hom  ( \mathbf{P}, F)
\rightarrow \Hom  (\mathbf{P}, Q')\rightarrow 0$$ is an exact
sequence of complexes (the exactness of $\Hom  ( \mathbf{P}, F)
\rightarrow \Hom  (\mathbf{P}, Q')\rightarrow 0$ follows by
\cite[Proposition 2.3]{HH}, since  $M$ is $n$-strongly Gorenstein
projective then Gorenstein projective). Now, by inductive
hypothesis, the complexes  $\Hom ( \mathbf{P}, K)$ and $ \Hom  (
\mathbf{P}, F)$ are exact, then so is $\Hom ( \mathbf{P}, Q')$ (by
\cite[Theorem 6.3]{Rot}), as desired.\\
To prove the equivalences  $1\Leftrightarrow 3 $ and
$2\Leftrightarrow 4 $, it suffices to apply the following remark. If
there exists an exact sequence of the form: $$ 0 \rightarrow
M\rightarrow P_n\rightarrow\cdots\rightarrow P_1 \rightarrow
M\rightarrow 0 ,$$ where each $P_i$ is a projective module, we get
$\Ext^{i}(M,G)\cong\Ext^{n+i}(M,G)$ for all $i\geq 1$ and all
modules $G$ (from \cite[Theorem 9.4]{Rot}). This implies that
$\Ext^{i}(M,G)=0$ for  all $i\geq 1$  whenever $n$ successive terms
of $\Ext(M,G)$ vanish.\cqfd\bigskip

Recall that a ring $R$ is called Iwanaga-Gorenstein ring (or
simply a Gorenstein ring), if $R$ is both left and right
Noetherian and if $R$ has finite self-injective dimension on  both
the left and the right \cite[Definition 9.1.1]{Rel-hom}. Over
Gorenstein rings, $n$-strongly Gorenstein projective modules have
the following characterization:

\begin{cor}\label{cor-SG-Gor-ring}
If $R$ is a Gorenstein ring, then an $R$-module $M$ is
$n$-strongly Gorenstein projective if  and only if there exits an
exact sequence of the form $0\rightarrow M\rightarrow
P_n\rightarrow\cdots\rightarrow P_1 \rightarrow M\rightarrow 0,$
where each $P_i$ is a projective $R$-module.
\end{cor}
\proof Simply apply Theorem \ref{thm-car-n-SG-proj} and the fact
that over a Gorenstein ring every projective module has finite
injective dimension \cite[Proposition 9.1.7]{Rel-hom}.\cqfd\bigskip

Now, we give change of rings results.

\begin{thm}\label{thm-ch-ring-SGp}
Let $R\rightarrow S$ be a homomorphism of commutative rings  with
$\pd_R(S)<\infty$. If $M$ is an $n$-strongly Gorenstein projective
$R$-module ($n\geq 1$), then $ S\otimes_R M$ is an $n$-strongly
Gorenstein projective $S$-module.
\end{thm}
\proof Let $M$ be an $n$-strongly Gorenstein projective $R$-module.
Then, there exists, by Theorem \ref{thm-car-n-SG-proj} and its
proof, an exact sequence of $R$-modules $0\rightarrow M\rightarrow
P_n\rightarrow\cdots\rightarrow P_1 \rightarrow M\rightarrow 0$,
where each $P_i$ is projective, and $\Ext_{R}^{i}(M,Q')=0$  for any
$i\geq 1$ and for any
$R$-module $Q'$ with finite projective dimension.\\
As the end of the proof of Theorem \ref{thm-car-n-SG-proj},
$\Tor_{i}^{R}(S,M)=0$ for all $i\geq 1$  whenever $n$ successive
terms of $\Tor^{R}(S,M)$ vanish. This holds since $\pd_R(S)<\infty$,
which shows that the following sequence of $S$-modules:
$$0\rightarrow S\otimes_R  M\rightarrow S \otimes_R P_n
 \rightarrow\cdots\rightarrow S  \otimes_R P_1 \rightarrow S
\otimes_R M\rightarrow 0$$ is exact. Note that each $ S \otimes_R
P_i$ is a projective $S$-module.\\
On the other hand, for any  $S$-module $N$, $\Ext_{S}^{i}(
S\otimes_R M, N)=\Ext_{R}^{i}(M,N)$ or all $i\geq 1$ (by
\cite[Proposition 4.1.3]{CE}). If the $S$-module $N$ is projective,
then it has finite projective dimension as an $R$-module. Indeed,
$\pd_R(S)<\infty$ implies that $\pd_R(S^{(I)})<\infty$ for any index
set $I$, and so does $N$ being  a direct summand of a free
$S$-module. Then, $\Ext_{S}^{i}( S\otimes_R M,
N)=\Ext_{R}^{i}(M,N)=0$ by the observation above. Therefore, by
Theorem \ref{thm-car-n-SG-proj},  $S  \otimes_R M $ is an
$n$-strongly Gorenstein projective $S$-module.\cqfd

\begin{cor}\label{cor-ch-ring-SGp}
Let $R\rightarrow S$ be a homomorphism of commutative rings with
$\pd_R(S)<\infty$. If $M$ is a Gorenstein projective $R$-module,
then $S \otimes_R M$ is a  Gorenstein projective $S$-module.
\end{cor}
\proof Let $M$ be a Gorenstein projective $R$-module. Then, it is
a direct summand of a strongly Gorenstein projective $R$-module
$G$ (by \cite[Theorem 2.7]{BM}). Then, $S \otimes_R M$  is a
direct summand of the $S$-module $S \otimes_R G$ which is, from
Theorem \ref{thm-ch-ring-SGp}, strongly Gorenstein projective.
Therefore, from \cite[Theorem 2.7]{BM}, $S \otimes_R M$ is a
Gorenstein projective $S$-module, as desired.\cqfd\bigskip

Finally, it is convenient to mention that one could define and
study the $n$-strongly injective and flat modules as we have done
in the Gorenstein projective case, and so Proposition
\ref{Pro-relation}, Theorems \ref{thm-car-n-SG-proj} and
\ref{thm-ch-ring-SGp}, and Corollaries \ref{cor-SG-Gor-ring} and
\ref{cor-ch-ring-SGp} have  the Gorenstein injective and flat
counterparts.\bigskip

\noindent {\bf Acknowledgements.} The authors thank the referee
for his/her careful reading of this work.

\end{section}



\begin{thebibliography}{999}\addcontentsline{toc}{section}{\protect\numberline{}{Bibliography}}


\bibitem{A1} M. Auslander, Anneaux de Gorenstein et torsion en algèbre commutative, Secrétariat mathématique, Paris, 1967, Séminaire d'algèbre commutative dirigé par Pierre Samuel, 1966/67. Texte rédigé, d'après des exposés de Maurice Auslander, par Marquerite Mangeney, Christian Peskine et Lucien Szpiro, Ecole Normale Superieure de Jeunes Filles.
\bibitem{A2} M. Auslander and M. Bridger,  \textit{Stable module theory}, Memoirs. Amer. Math. Soc.  \textbf{94} (American Mathematical Society, Providence,  Rhode Island, 1969).
\bibitem{BM} D. Bennis and N. Mahdou,  Strongly Gorenstein projective, injective, and flat modules, \textit{J. Pure Appl. Algebra} \textbf{210} (2007) 437--445.
\bibitem{Period} D. J. Benson and K. R. Goodearl, Periodic flat modules, and flat modules for finite Groups, \textit{Pacific J. Math.} \textbf{196} (2000) 45--67.
\bibitem{CE}  H. Cartan and S. Eilenberg, \textit{Homological Algebra}, Princeton Mathematical Series \textbf{19} (Princeton University Press, Princeton, 1956).
\bibitem{LW} L. W. Christensen, \textit{Gorenstein dimensions}, Lecture Notes in Math. \textbf{1747} (Springer-Verlag, Berlin, Heidelberg, 2000).
\bibitem{CFH} L. W. Christensen, A. Frankild, and H. Holm, On Gorenstein projective, injective and flat dimensions - a functorial description with applications, \textit{J. Algebra} \textbf{302} (2006) 231--279.
\bibitem{Rel-hom} E. E. Enochs and O. M. G. Jenda, \textit{Relative homological algebra}, de Gruyter Expositions in Mathematics \textbf{30} (Walter de Gruyter \& Co., Berlin, 2000).
\bibitem{GoInPj}E. Enochs and O. Jenda, Gorenstein injective and projective modules, \textit{Math. Z.} \textbf{220} (1995) 611--633.
\bibitem{GoIn}  E. E. Enochs and  O. M. G. Jenda, On Gorenstein injective modules, \textit{Comm. Algebra} \textbf{21}  (1993) 3489--3501.
\bibitem{GoPlat} E. Enochs, O. Jenda and B. Torrecillas, Gorenstein flat modules, \textit{Nanjing Daxue Xuebao Shuxue Bannian Kan} \textbf{10} (1993) 1--9.
\bibitem{Fox}  E. Enochs, O. Jenda and J. Xu, Foxby duality and Gorenstein injective and projective modules, \textit{Trans. Amer. Math. Soc.} \textbf{348} (1996) 3223--3234.
\bibitem{HH} H. Holm, Gorenstein homological dimensions, \textit{J. Pure Appl. Algebra} \textbf{189} (2004)  167--193.
\bibitem{Rot}       J. J. Rotman, \textit{An Introduction to Homological Algebra} (Academic Press, London, New York, 1979).
\end{thebibliography}
\end{document}